\documentclass[12pt,titlepage]{article}
\usepackage[applemac]{inputenc}
\usepackage{anysize}
\usepackage{float}
\usepackage{amsmath}
\usepackage{amssymb}
\usepackage{amsbsy}
\usepackage{amsfonts}
\usepackage{amsthm}
\usepackage{latexsym}

\def\Q{\Bbb Q}
\def\Z{\Bbb Z}
\def\F{\Bbb F}
\def\AA{\mathcal A}
\def\DD{\mathcal D}
\def\QQ{\mathcal Q}

\def\Aut{\mathop{\rm Aut}\nolimits}

\def\Aff{\mathop{\rm Aff}\nolimits}

\def\ker#1{{\rm ker}(#1)}

\let\phi=\varphi

\newcommand{\apl}[3]{#1\colon #2\longrightarrow #3}

\newcommand{\Irr}{\operatorname{Irr}}

\newtheorem{lemma}{Lemma}[section]
\newtheorem{proposition}[lemma]{Proposition}
\newtheorem{theorem}[lemma]{Theorem}

\newtheorem{definition}{Definition}[section]

\begin{document}

~\bigskip

\thispagestyle{empty}

 \centerline{\Large Words and characters in finite $p$-groups}

\bigskip
\centerline{by}
\bigskip

\centerline{\bf Ainhoa I\~niguez}
\centerline{Mathematical Institute} \centerline{University of Oxford}
 \centerline{Andrew Wiles Building}
  \centerline{Woodstock Road}
   \centerline{OX2 6GG, Oxford}
    \centerline{UNITED KINGDOM}
 \centerline{E-mail: ainhoa.iniguez@maths.ox.ac.uk}
\smallskip
\centerline{and}
\smallskip
\centerline{\bf Josu Sangroniz\footnote{Corresponding
author}}\centerline{Departamento de Matem\'aticas}
\centerline{Facultad de Ciencia y Tecnolog\'\i a}
\centerline{Universidad del Pa\'{\i}s Vasco UPV/EHU} \centerline{48080
Bilbao} \centerline{SPAIN} \centerline{E-mail:
josu.sangroniz@ehu.es}

\bigskip
\bigskip
\centerline{\bf Abstract}
\bigskip
 \centerline{\parbox{13cm}{
 Given a group word $w$ in $k$ variables, a finite group $G$ and $g\in G$, we consider the number $N_{w,G}(g)$ of $k$-tuples $g_1,\dots ,g_k$ of elements of $G$ such that $w(g_1,\dots ,g_k)=g$. In this work we study the  functions $N_{w,G}$ for the class of nilpotent groups of nilpotency class $2$. We show that, for the groups in this class, $N_{w,G}(1)\geq |G|^{k-1}$, an inequality that can be improved to $N_{w,G}(1)\geq |G|^k/|G_w|$ ($G_w$ is the set of values taken by $w$  on $G$) if $G$ has odd order. This last result is explained by the fact that the functions $N_{w,G}$ are  characters of $G$ in this case. For groups of even order, all that can be said is that $N_{w,G}$ is a generalized character, something that is false in general for  groups of nilpotency class greater than $2$. We characterize group theoretically when $N_{x^n,G}$ is a character if $G$ is a $2$-group of nilpotency class $2$. Finally we also address the (much harder) problem of studying if $N_{w,G}(g)\geq |G|^{
 k-1}$ for $g\in G_w$, proving that this is the case for the free $p$-groups of nilpotency class $2$ and exponent $p$.}}

\bigskip
\centerline{\parbox{13cm}{\noindent Keywords: $p$-groups; words; characters}}
\bigskip
\centerline{\parbox{13cm}{\noindent MSC: 20D15, 20F10}}

\bigskip\bigskip

\noindent{\footnotesize Both authors are supported by the MINECO (grants MTM2011-28229-C02-01 and MTM2014-53810-C2-2-P).  The second author is also supported by the Basque Government (grants IT753-13 and IT974-16).}

\newpage
\section{Introduction}

Given a group word $w$ in $k$ variables $x_1,\dots ,x_k$, that is an element in the free group $F_k$ on $x_1,\dots ,x_k$,   we can define for any $k$ elements $g_1,\dots ,g_k$ in a group $G$ the element $w(g_1,\dots ,g_k)\in G$ by applying to $w$ the group homomorphism from $F_k$ to $G$ sending $x_i$ to $g_i$. For $G$ a finite group and $g\in G$ we consider the number
\begin{equation}\label{def:N}N_{w,G}(g)=|\{(g_1,\dots ,g_k)\in G^{(k)}\mid w(g_1,\dots ,g_k)=g\}|.\end{equation}
($G^{(k)}$ denotes the $k$-fold cartesian product of $G$ with itself.) So $N_{w,G}(g)$ can be thought of as the number of solutions of the equation $w=g$. The set of word values of $w$ on $G$, i.~e.,~the set of elements $g\in G$ such that the equation $w=g$ has a solution in $G^{(k)}$, is denoted $G_w$.

There is an extensive literature on the functions $N_{w,G}$, sometimes expressed in terms of the probabilistic distribution $P_{w,G}=N_{w,G}/|G|^k$. For example Nikolov and Segal gave in \cite{NikolovSegal} a characterization of the finite nilpotent (and also solvable) groups based on   these probabilities: $G$ is nilpotent if and only if $\inf_{w,g}P_{w,G}(g)>p^{-|G|}$, where $w$ and $g$ range over all words and $G_w$, respectively, and $p$ is the largest prime dividing $|G|$. One of the implications is easy: if $G$ is not nilpotent the infimum is in fact zero. Indeed, we can consider the $k$-th left-normed lower central word $\gamma_k=[[[x_1,x_2],x_3],...,x_k]$. Since $G$ is not nilpotent, there exists some non-trivial element $g\in G_{\gamma_k}$ (for any $k$). Since $\gamma_k(g_1,\dots,g_k)=1$ if some $g_i=1$, we have that $P_{\gamma_k,G}(g)\leq (|G|-1)^k/|G|^k$, which can be made arbitrarily small. On the other hand the estimation $P_{w,G}(g)>p^{-|G|}$ for $g\in G_w$ seems to be 
 far from sharp and Amit in an unpublished work has conjectured that if $G$ is   nilpotent, $P_{w,G}(1)\geq 1/|G|$.

We prefer to give our results in terms of the functions $N_{w,G}$. In this paper we focus our attention on finite nilpotent groups of nilpotency class $2$, which we take to be $p$-groups right from the outset, so all the results quoted here are referred to this family of groups. In Section 2 we consider a natural equivalence relation for words that enable us to assume that they have a special simple form. Then it is not difficult to prove Amit's conjecture $N_{w,G}(1)\geq |G|^{k-1}$ for $w\in F_k$. This result has been proved independently by Levy in \cite{Levy} using a similar procedure, although our approach to the concept of word equivalence is different. In the next two sections we show that the functions $N_{w,G}$ are generalized characters, a result that is false for nilpotent groups of nilpotency class greater than $2$,  and even more, if $G$ has odd order, they are genuine characters. In particular we obtain an improvement of Amit's conjectured bound, namely, $N_{w,G}(1)\geq 
 |G|^k/|G_w|$. For $2$-groups, there are easy examples where  $N_{x^2,G}$ fails to be a character and we actually characterize group-theoretically when this happens for the power words $w=x^n$ (always within the class of nilpotent $2$-groups of nilpotency class $2$). In the last section we consider briefly the conjecture  $N_{w,G}(g)\geq |G|^{k-1}$ for $w\in F_k$ and $g\in G_w$. This problem is much harder than the case $g=1$ and only some partial results have been obtained, for instance confirming the conjecture if $G$ is a free nilpotent $p$-group of nilpotency class $2$ and exponent $p$.

\section{Words in $p$-groups of nilpotency class $2$}

Since we are going to work with $p$-groups of nilpotency class $2$, it is more convenient for us to define a \emph{word} in the variables $x_1,\dots, x_k$ as an element in the free pro-$p$ group of nilpotency class $2$ on the variables $x_1,\dots ,x_k$, $F_k$. Thus, if $w\in F_k$ is a word, it can be represented in a unique way as
$$w=x_1^{z_1}\dots x_k^{z_k}\prod_{1\leq i<j\leq k}[x_i,x_j]^{z_{ij}},$$
where the exponents $z_l$, $z_{ij}$ are $p$-adic integers. Of course, if $G$ is a finite $p$-group (or pro-$p$ group) of nilpotency class $2$ and $g_1,\dots,g_k\in G$, it makes sense to evaluate $w$ on $g_1,\dots ,g_k$  by   applying the homomorphism $\apl \pi{F_k}G$ given by $x_i\mapsto g_i$. As in the introduction, we denote this element $w(g_1,\dots,g_k)$ and   define the function $N_w=N_{w,G}$ by (\ref{def:N}).

If $\sigma$ is an automorphism of $F_k$, $\sigma$ is determined by the images of the generators $x_1,\dots ,x_k$, which we denote $w_1,\dots,w_k$. Then the image of $w\in F_k$ is the word $w(w_1,\dots,w_k)$, the result of evaluating $w$ on $w_1,\dots,w_k$. Since $\sigma$ is an automorphism, there exist $w_1',\dots ,w_k'\in F_k$ such that $w_i'(w_1,\dots,w_k)=x_i$, $1\leq i\leq k$, and the inverse automorphism is given by $x_i\mapsto w_i'$. If $G$ is a finite $p$-group (or pro-$p$ group) of nilpotency class $2$, we can define the map $\apl\varphi{G^{(k)}}{G^{(k)}}$ by $\varphi (g_1,\dots ,g_k)=(w_1(g_1,\dots ,g_k),\dots ,w_k(g_1,\dots ,g_k))$ and it is clear that this map is a bijection with the inverse map given by $(g_1,\dots,g_k)\mapsto (w_1'(g_1,\dots ,g_k),\dots ,w_k'(g_1,\dots ,g_k))$. If $w'=w(w_1,\dots,w_k)$, it is clear that $w'(g_1,\dots,g_k)=g$ if and only if $w(\varphi (g_1,\dots ,g_k))=g$, thus $\phi$ is a bijection between the solutions of $w'=g$ and $w=g$ and  in partic
 ular, $N_{w,G}=N_{w',G}$.

\begin{definition} We will say that two words $w,w'\in F_k$ are \emph{equivalent} if they belong to the same orbit under the action of the automorphism group of $F_k$.
\end{definition}

Therefore we have proved the following result.

\begin{proposition} If $w,w'\in F_k$ are two equivalent words, $N_{w,G}=N_{w',G}$ for any finite $p$-group $G$ of nilpotency class $2$.\end{proposition}

Next we want to find a set of representatives of the equivalence classes of words. There are natural homomorphisms
\begin{equation}\label{homs}\Aut (F_k)\twoheadrightarrow\Aut (F_k/F_k')\cong GL_k(\Z_p)\rightarrow\Aut(F_k'),\end{equation}
where the composite map is the restriction. For the middle isomorphism, given an automorphism,
 we write the coordinates of the images of the vectors in a basis of the $\Z_p$-module $F_k/F_k'\cong \Z_p^{(k)}$ as rows of the corresponding matrix. The last homomorphism  comes from the identification of $F_k'$ with the exterior square of $F_k/F_k'$. More explicitly, we identify $F_k'$ with the additive group of the $k\times k$ antisymmetric matrices over $\Z_p$, $\AA_k$, via $w=\prod_{i<j}[x_i,x_j]^{z_{ij}}\mapsto A$, where $A\in \AA_k$ has entries $z_{ij}$ for $1\leq i<j\leq k$. Then, if $X\in GL_k(\Z_p)$, the action of $X$ on $\AA_k$ is given by $A\mapsto X^tAX$. This action is better understood if we interpret $A$ as an alternating bilinear form on the free $\Z_p$-module $\Z_p^{(k)}$. Notice however that, under a change of basis, the matrix  $A$ is now transformed into $PAP^t$, where $P$ is the  matrix associated to the change of basis, writing the coordinates of the vectors in the new basis as rows of $P$.

 We start by analyzing the action of $GL_k(\Z_p)$ and the affine subgroup
 $$\Aff_{k-1}(\Z_p)=\left\{
\begin{pmatrix}1&0\\u^t&X\end{pmatrix} \mid   u\in\Z_p^{(k-1)},  X\in GL_{k-1}(\Z_p)\right\}$$
 ($t$ means transposition), on $\AA_k$, giving a set of representatives of the orbits. The result about the action of $GL_k(\Z_p)$ generalizes naturally if we  replace $\Z_p$ by any principal ideal domain (see for example \cite[Theorem IV.1]{Newman}), but a more elaborate   proof is required. For completeness we include a proof here  that  takes advantage of the fact that $\Z_p$ is a discrete valuation ring and can be later adapted to the case when the acting group is $\Aff_{k-1}(\Z_p)$.

\begin{lemma}\label{canonicalmatrix}
\begin{enumerate}
\item Each orbit of the action of  $GL_k(\Z_p)$ on $\AA_k$ contains a unique  block diagonal   matrix with diagonal non-zero blocks $p^{s_i}H$, $H=\begin{pmatrix}\phantom{-}0&1\\-1&0\end{pmatrix}$, $1\leq i\leq r$  and  $0\leq s_1\leq\dots\leq s_r$  {\rm (}$0\leq r\leq k/2${\rm )}.
\item  Each orbit of the action  of the affine group $\Aff_{k-1}(\Z_p)$ on $\AA_k$ contains a unique  tridiagonal matrix $A$ {\rm (}that is, all the entries $a_{ij}$ of $A$ with $|i-j|>1$ are zero{\rm )} with  the non-zero entries above the main diagonal $a_{i,i+1}=p^{s_i}$, $1\leq i\leq r$, and $0\leq s_1\leq s_2\leq \dots \leq s_r$ {\rm (}$0\leq r< k${\rm )}.
\end{enumerate}
\end{lemma}

\begin{proof}
Given $A$ in $\AA_k$, we consider a basis $\{e_1,\dots ,e_k\}$ (for instance, the canonical basis) in the free $\Z_p$-module $\Z_p^{(k)}$ and the alternating bilinear form $(\  ,\, )$ defined by the matrix $A$ with respect to this basis. There is nothing to prove if $A$ is the zero matrix, so we can suppose that $(e_i,e_j)\not=0$ for some $1\leq i<j\leq k$ and we can assume that its $p$-adic valuation  is   minimum among the valuations of all the (non-zero) $(e_r,e_s)$. After reordering the basis, we can suppose that $i=1$ and $j=2$ and moreover, by multiplying $e_1$ or $e_2$ by a $p$-adic unit, we can suppose that $(e_1,e_2)=p^{s_1}$ for some $s_1\geq 0$. Notice that any (non-zero) $(u,v)$ has $p$-adic valuation greater than or equal to $p^{s_1}$.

Now for each $i\geq 3$ we   set $e_i'=e_i+\alpha_i e_1+\beta_i e_2$, where $\alpha_i,\beta_i\in\Z_p$ are chosen so that $(e_i',e_1)=(e_i',e_2)=0$. The elements $\alpha_i,\beta_i$ exist because the valuation of $(e_1,e_2)$ is less than or equal to the valuations of $(e_i,e_2)$ and $(e_i,e_1)$. By replacing $e_i$ by $e_i'$ we can suppose that $\langle e_1,e_2\rangle$ is orthogonal to $\langle e_3,\dots ,e_k\rangle$. Proceeding inductively, we obtain a basis $\{e'_1,\dots,e'_k\}$ such that, for some $0\leq r\leq k/2$, the subspaces $\langle e'_{2i-1},e'_{2i}\rangle$ are pairwise orthogonal for $1\leq i\leq r$, the remaining vectors are in the kernel of the form and  $(e'_{2i-1},e'_{2i})=p^{s_i}$, $1\leq i\leq r$, with $0\leq s_1\leq\dots\leq s_r$.
It is clear that with respect to this new basis the matrix associated to the form $(\  ,\, )$   has the desired form.

To prove uniqueness   suppose that $A$ and $A'$ are block diagonal  matrices with (non-zero) diagonal blocks $p^{s_1}H,\dots ,p^{s_r}H$, $0\leq s_1\leq\dots\leq s_r$, and $p^{s'_1}H,\dots ,p^{s'_t}H$, $0\leq s'_1\leq\dots\leq s'_t$, respectively, and $A'=X^tAX$ for some $X\in GL_k(\Z_p)$. The matrices $A$, $A'$ and $X$ can be viewed as endomorphisms of the abelian groups $R_n=(\Z/p^n\Z)^{(k)}$, $n\geq 1$. Since   $X$   defines in fact an automorphism of $R_n$ the image subgroups of $A$ and $A'$ (as endomorphisms of $R_n$) have the same order. For $A$ this order is $p^{2s}$, where $s=\sum_{s_i\leq n}(n-s_i)$, and similarly for $A'$. We conclude that, for any $n\geq 1$, $\sum_{s_i\leq n}(n-s_i)=\sum_{s'_i\leq n}(n-s'_i)$, whence $r=t$ and $s_i=s'_i$ for all $1\leq i\leq r$, that is, $A=A'$.

For the existence  part in (ii) we have to show that, given an alternating form $(\  ,\, )$ on $\Z_p^{(k)}$ and a basis $\{e_1,\dots ,e_k\}$, there exists another basis $\{e'_1,\dots ,e'_k\}$ such that $e'_1\in e_1+\langle e_2,\dots ,e_k\rangle$, $\langle e'_2,\dots ,e'_k\rangle=\langle e_2,\dots ,e_k\rangle$ and $(e'_i,e'_j)=0$ for $|i-j|>1$, $(e_i,e_{i+1})=p^{s_i}$, $0\leq s_1\leq \dots \leq s_r$,  $(e_i,e_{i+1})=0$, $r<i<k$. We can suppose that $(\  ,\, )$ is not the trivial form and then consider the minimum valuation $s_1$ of all the (non-zero) $(e_i,e_j)$. If this minimum is attained for some $(e_1,e_j)$ we interchange $e_2$ and $e_j$. Otherwise this minimum is attained for some $(e_i,e_j)$, $2\leq i<j\leq k$ and $(e_1+e_i,e_j)$ has still valuation $s_1$ (because the valuation of $(e_1,e_j)$ is strictly greater than $s_1$). By replacing $e_1$ by $e_1+e_i$, interchanging $e_2$ and $e_j$, and adjusting  units, we can suppose that $(e_1,e_2)=p^{s_1}$. Now we can replace $e_i$, $i\
 geq 3$, by $e'_i=e_i+\alpha_ie_2$, where $\alpha_i$ is chosen so that $(e_1,e'_i)=0$. Thus we can assume $(e_1,e_i)=0$ for $i\geq 3$. Now we iterate the procedure with the basis elements $e_2,\dots ,e_k$.

We prove uniqueness with a similar counting argument as in (i) but this time considering the order of the images of the subgroup of $R_n$, $S_n=\{0\}\times (\Z/p^n\Z)^{(k-1)}$. So we assume that $A'=X^tAX$ with $A$ and $A'$ as in (ii) and $X\in\Aff_{k-1}(\Z_p)$. Notice that, as an automorphism of $R_n$, $X$ fixes $S_n$, so the images of $S_n$ by $A$ and $A'$ must have the same order. These orders are $p^s$ and $p^{s'}$, where $s=\sum_{s_i\leq n}(n-s_i)$ and similarly for $s'$, so $s=s'$ and, since this must happen for any $n\geq 1$, it follows that $s_i=s'_i$ for all $i$, that is, $A=A'$.\end{proof}

\begin{proposition}\label{prp:canonicalword} The following words are a system of representatives of the action of $\Aut F_k$ on $F_k$:
\begin{eqnarray}
&&\label{commcanonicalword}[x_1,x_2]^{p^{s_1}}\dots [x_{2r-1},x_{2r}]^{p^{s_r}},\ 0\leq r\leq k/2,\  0\leq s_1\leq\dots\leq s_r,\\
&&\label{noncommcanonicalword}x_1^{p^{s_1}}[x_1,x_2]^{p^{s_2}}[x_2,x_3]^{p^{s_3}}\dots [x_{r-1},x_r]^{p^{s_r}},\ 1\leq r\leq k,\  s_1\geq 0,\  0\leq s_2\leq\dots\leq s_r.
\end{eqnarray}
\end{proposition}

\begin{proof}
As explained above, the action of $\Aut F_k$ on $F_k'$ can be suitably identified with the action of $GL_k(\Z_p)$ on $\AA_k$, thus it follows directly from the part (i) of the previous lemma that the words  (\ref{commcanonicalword}) are representatives for the orbits contained in $F_k'$.

Now suppose $w\in F_k\backslash F_k'$. Then
$w=(x_1^{z_1}\dots x_k^{z_k})^{p^{s_1}}\prod_{1\leq i<j\leq k}[x_i,x_j]^{z_{ij}}$, where $s_1\geq 0$ and some $z_i$ is a $p$-adic unit. After applying the inverse of the automorphism $x_1\mapsto x_1^{z_1}\dots x_k^{z_k}$, $x_i\mapsto x_1$, $x_j\mapsto x_j$, for $j\not=1,i$, we can assume    that $x_1^{z_1}\dots x_k^{z_k}=x_1$.  Now we consider the action of the stabilizer of $x_1$, $\Aut_{x_1}F_k$. The image of this subgroup by the first map in (\ref{homs}) is $\Aff_{k-1}(\Z_p)$, so we can identify the action of $\Aut_{x_1}F_k$ on $F_k'$ with the action of $\Aff_{k-1}(\Z_p)$ on $\AA_k$. It follows from Lemma \ref{canonicalmatrix} (ii) that $w$ is equivalent to a word  in (\ref{noncommcanonicalword}). Notice also that if $w'=\sigma(w)$ for two of these words $w$ and $w'$ and some $\sigma\in\Aut F_k$, it would follow by passing to $F_k/F_k'$ that $\sigma(\overline x_1)^{p^{s_1}}=\overline x_1^{p^{s'_1}}$. Since $\sigma$ induces automorphisms of $(F_k/F_k')^{p^s}$ and this chain of subg
 roups of $F_k/F_k'$ is strictly decreasing, we conclude that $s_1=s'_1$. But $F_k/F_k'$ is torsion-free, so $\sigma$ fixes $\overline x_1$, that is, $\sigma(x_1)=x_1z$ for some $z\in F_k'$. Composing $\sigma$ with the automorphism $x_1\mapsto x_1z^{-1}$, $x_i\mapsto x_i$, $i\geq 2$, we can suppose that $\sigma\in \Aut_{x_1}F_k$. Thus, the two matrices in $\AA_k$ associated to $x_1^{-p^{s_1}}w$ and $x_1^{-p^{s_1}}w'$ are in the same orbit by $\Aff_{k-1}(\Z_p)$, and so they coincide by Lemma \ref{canonicalmatrix} (ii). We conclude that $w=w'$.
\end{proof}

\begin{theorem} Let $G$ be  a finite $p$-group of nilpotency class $2$ and $w$ a word in $k$ variables. Then $N_w(1)\geq |G|^{k-1}$.
\end{theorem}

\begin{proof}
We can suppose that $w$ is as in the last proposition. Write $k_0=\lfloor k/2\rfloor$ and fix $g_2,g_4\dots ,g_{2k_0}\in G$. Then the map $G'\times G^{(k-k_0-1)}\longrightarrow G'$ given by $(x_1,x_3,\dots ,x_{2(k-k_0)-1})\mapsto w(x_1,g_2,x_3,\dots)$ is a group homomorphism whose kernel has size at least $|G|^{k-k_0-1}$. Since there are $|G|^{k_0}$ choices for $g_2,g_4,\dots ,g_{2k_0}$, we get at least $|G|^{k-1}$ solutions to the equation $w(x_1,\dots,x_k)=1$.
\end{proof}

\section{The functions $N_w$ from a character-theoretical point of view}

In this section, unless otherwise stated, we consider an arbitrary
finite group $G$ and a word $w$ that is thought  now as an element
in the free group with, say, free generators $x_1,\dots,x_k$. The
function $N_w=N_{w,G}$ is of course a class function, so it can be
written as a linear combination of the irreducible characters of
$G$:
\begin{equation}\label{decomposition}N_w=\sum_{\chi\in\Irr (G)}N_w^\chi\chi,\end{equation} where
\begin{equation}\label{coefficients}N_w^\chi=(N_w,\chi)=\frac{1}{|G|}\sum_{g\in
G}N_w(g)\overline{\chi(g)}=\frac{1}{|G|}\sum_{(g_1,\dots ,g_k)\in
G^{(k)}}\overline{\chi(w(g_1,\dots,g_k))}.\end{equation}

It is a natural question to study when the function $N_w$ is a
character  of $G$. Notice that in this case $N_w(g)\leq N_w(1)$ for
any element $g\in G$, so $N_w(1)$ will  be at least the average of
the non-zero values of the function $N_w$, that is
\begin{equation}\label{improvedineq}N_w(1)\geq \frac{1}{|G_w|}\sum_{g\in G_w}N_w(g)=\frac{|G|^k}{|G_w|},\end{equation}
which of course improves the bound conjectured by Amit, $N_w(1)\geq |G|^{k-1}$.

It is easy to find examples where $N_w$ is not a character. Probably
the simplest is $Q_8$, the quaternion group of order $8$, and
$w=x^2$: $N_{x^2,Q_8}(1)=2$ but $N_{x^2,Q_8}(z)=6$ for the unique
involution $z\in Q_8$. For $p$ odd we can construct a $p$ analogue to $Q_8$ as
$G=\langle g_1,\dots,g_{p-1}\rangle\rtimes\langle
g_{p}\rangle/\langle g_{p-1}^pg_{p}^{-p}\rangle$, where  $\langle
g_1,\dots,g_{p-1}\rangle\cong C_p\times\dots\times C_p\times
C_{p^2}$, $\langle g_{p}\rangle\cong C_{p^2}$  ($C_n$ denotes a cyclic group of order $n$) and
$g_i^{g_{p}}=g_ig_{i+1}$, $1\leq i<p-2$,
$g_{p-2}^{g_{p}}=g_{p-2}g_{p-1}^p$ and
$g_{p-1}^{g_p}=g_{p-1}g_1^{-1}$. It can be checked that
$N_{x^p,G}(1)=p^{p-1}$ but  $N_{x^p,G}(z)=p^p+p^{p-1}$ for any
non-trivial   element $z\in Z(G)=G_{x^p}=\langle g_p^p\rangle$. Notice that $|G|=p^{p+1}$ and this
is the smallest order for which $N_{x^p,G}$ can fail to be a
character, since $p$-groups of order at most $p^p$ are regular and,
for regular $p$-groups, $N_{x^p,G}$ is the regular character of
$G/G^p$. Also, for the quaternion group and its $p$ analogues, (\ref{improvedineq}) is false.
However, it is a well known result that  in general  the
functions  $N_{x^n,G}$ are generalized characters
(i.~e.,~$\Z$-linear combinations of irreducible characters), see
\cite[Problem 4.7]{Isaacs}.

For words $w$ in more than one variable the situation is different
and there are examples of groups  $G$ and words $w$ where $N_{w,G}$
is not a generalized character, even among  nilpotent groups. As for non-solvable examples one can take $G=PSL_2(11)$
and the $2$-Engel word $w=[x,y,y]$ (see \cite{Parzanchevski} for another choice of $w$). Then the coefficients $N_w^\chi$ for the
two irreducible characters $\chi$ of degree $12$ are $305\pm 23\sqrt 5$.
More examples can be obtained using the following result by A.
Lubotzky \cite{Lubotzky}: if $G$ is a simple group and $1\in
A\subseteq G$ is a subset invariant under the group of automorphisms
of $G$, then $A=G_w$ for some word $w$ in two variables. Notice that
if $A$ contains an element $a$ such that $a^i\not\in A$, for some
$i$ coprime with the order of $a$, then $N_w(a^i)=0$ but
$N_w(a)\not=0$, something that cannot happen if $N_w$ is a
generalized character.

Examples of $p$-groups where
$N_w$ is not a generalized character are provided by the free $p$-groups of nilpotency class $4$ and exponent $p$ with $p>2$ and $p\equiv 1\pmod 4$ and $w=[x,y,x,y]$ (which settles in the negative a question of Parzanchevski who had asked if the functions $N_w$ were always generalized characters for solvable or nilpotent groups). We realize these groups as as $1+J$, where $J=I/I^4$ and $I$ is the ideal
generated by $X$ and $Y$ in the algebra of the polynomials in the
non-commuting indeterminates $X$ and $Y$ with coefficients in
$\F_p$. If $x=1+X$ and $y=1+Y$ and $u=[x,y,x,y]$, then certainly $u\in G_w$ but we claim that if $i$ is not a quadratic residue modulo $p$, then $u^i\not\in G_w$ (so $N_w(u)\not=0$ but $N_w(u^i)=0$). Indeed, one can check directly (or by
appealing to the Lazard correspondence, but in this case, only for $p>3$, see \cite[\S 10.2]{Khukhro}) that
$\gamma_4(G)=1+\gamma_4(J)$, where $\gamma_4(J)$ is the fourth term
in the descending central series of the Lie algebra $J$. Now $u^i\in
G_w$ if and only if
$$i[X,Y,X,Y]=[aX+bY,cX+dY,aX+bY,cX+dY]$$
for some $a,b,c,d\in \F_p$, and one can see that this equation has no
solution if $i$ is not a quadratic residue modulo $p$.

In contrast to the previous example, we shall show  that the
functions $N_w$ are always generalized characters for $p$-groups of
nilpotency class $2$ and, in the next section, that they are
actually genuine characters for odd $p$. Before that we recall
briefly that for some words $w$ the functions $N_w$ are known to be
characters for any group, notably, for the commutator word $w=[x,y]$
(this is basically \cite[Problem 3.10]{Isaacs}). This classical result due to Frobenius can be
extended in various ways: when $w$ is an admissible word
(i.~e.,~a word in which all the variables appear exactly twice, once with
exponent $1$ and once with $-1$) \cite{DasNath} or when $w=[w',y]$, where $y$
is a variable which does not occur in $w'$ (so in particular, for
$\gamma_k=[x_1,x_2,\dots ,x_k]$, the $k$-th left-normed lower central word)
\cite{Strunkov}. It is also clear that if $N_w$ and $N_{w'}$ are
characters (or generalized characters), so is $N_{ww'}$ if $w$ and
$w'$ have no common variables. The reason is of course that
$N_{ww'}=N_w*N_{w'}$ is the convolution of the two functions $N_w$
and $N_{w'}$. More   information in this direction is given in \cite{Parzanchevski}.

As promised we finish this section by proving that the functions
$N_{w}$ are generalized characters for $p$-groups of nilpotency
class $2$. We   give first a  characterization  of when
$N_w$ is a generalized character. We have already used before that a
necessary condition is that $N_w(g)=N_w(g^i)$ for any $i$ coprime
with the order of $g$ and we are going to see that this condition is in fact
sufficient. The proof is standard once we know that the coefficients
$N_w^\chi$ in (\ref{decomposition}) are algebraic integers (as Amit
and Vishne point out in \cite{AmitVishne} this follows immediately
from (\ref{coefficients}) and the result of Solomon's in \cite{Solomon} that $N_w(g)$ is
always a multiple of $|C_G(g)|$).

\begin{lemma}\label{lem:coprime} Let $G$ be a group and $w$ a word. Then $N_w=N_{w,G}$ is
a generalized character of $G$ if and only if $N_w(g)=N_w(g^i)$ for
any $g\in G$ and $i$ coprime with the order of $G$.\end{lemma}

\begin{proof}
We only need to prove sufficiency. Since the coefficients $N_w^\chi$
are algebraic integers it suffices to see that they are rational
numbers. But, by elementary character and Galois theory, if $f$ is a
rational-valued class function of  a group $G$, $f$ is a $\Q$-linear combination of
irreducible characters if and only if $f(g)=f(g^i)$ for any $g\in G$
and $i$ coprime with the order of $G$. Indeed, if $f=\sum_{\chi\in\Irr (G)}a_\chi\chi$ with $a_\chi\in\Q$, and $\sigma$ is the automorphism of the cyclotomic extension $\Q(\varepsilon)/\Q$ sending $\varepsilon$ to $\varepsilon^i$, where $\varepsilon$ is a primitive $|G|$-th root of unity, we have
$$f(g)=f(g)^\sigma=\sum_{\chi\in\Irr(G)}a_\chi\chi^\sigma(g)=\sum_{\chi\in\Irr(G)}a_\chi\chi (g^i)=f(g^i).$$
And conversely, if $f(g)=f(g^i)$,
$$f(g)=f(g)^{\sigma^{-1}}=f(g^i)^{\sigma^{-1}}=(\sum_{\chi\in\Irr(G)}a_\chi\chi (g^i))^{\sigma^{-1}}=
\sum_{\chi\in\Irr(G)}a_\chi^{\sigma^{-1}}\chi (g).$$
By the linear independence of the irreducible characters, we conclude that $a_\chi=a_\chi^{\sigma^{-1}}$ for any automorphism $\sigma$, so $a_\chi\in\Q$.
\end{proof}

\begin{theorem}\label{thm:generalizedcharacter} Let $G$ be a $p$-group of nilpotency class $2$ and
$w$ a word. Then the function $N_w=N_{w,G}$ is a generalized
character of $G$.
\end{theorem}

\begin{proof}
By Proposition \ref{prp:canonicalword}   we can suppose that $w$ has the
form (\ref{commcanonicalword}) or (\ref{noncommcanonicalword}). Now we observe that, if $i$ is not a
multiple of $p$, the map $(g_1,g_2,\dots,g_k)\mapsto
(g_1^i,g_2,g_3^i,\dots)$ is a bijection from the set of solutions of
$w=g$ to the set of solutions of $w=g^i$, so in particular
$N_w(g)=N_w(g^i)$ and the result follows from the previous lemma.
\end{proof}

\section{The functions $N_w$ for odd $p$-groups of nilpotency class $2$}

The goal of this section is to show that $N_{w,G}$ is a genuine
character of a $p$-group $G$ of nilpotency class 2 if $p$ is odd. We
begin with a general result.

\begin{lemma}\label{coefswithkernel} Let $\chi\in\Irr (G)$ with kernel $K$ and $w$ a word in $k$ variables. Then $N_{w}^\chi=|K|^{k-1}N_{w,G/K}^{\overline\chi}$, where $\overline\chi$ is the character of $G/K$ defined naturally by $\chi$.
\end{lemma}

\begin{proof}
We have
$$N_{w}^\chi=(N_w,\chi)=(\tilde N_w,\overline\chi)_{G/K},$$
where $\tilde N_w$ is the average function defined  by $\tilde N_w
(g)=\frac{1}{|K|}\sum_{n\in K} N_{w}(gn)$ viewed as a function  on
$G/K$. As such a function it is clear that $\tilde
N_w=|K|^{k-1}N_{w,G/K}$, so the result is clear.
\end{proof}

We assume now that $G$ is a $p$-group of nilpotency class $2$. The
following technical result characterizes when $N_w$ is a character.

\begin{lemma}\label{basic} Let $G$ be a $p$-group of nilpotency class $2$.  Then $N_w$ is a character of $G$ if and only if for any (non-trivial) epimorphic image of $G$, say $G_1$, with cyclic center, $N_{w,G_1}(1)\geq N_{w,G_1}(z)$, where $z$ is a central element of $G_1$ of order $p$.
\end{lemma}

\begin{proof}
By Theorem \ref{thm:generalizedcharacter} we know that $N_w$ is a
generalized character, that is, all the numbers $N_w^\chi$ are
integers, so  $N_w$ is a character of $G$ if and only if these
numbers are non-negative.

We recall that a group $G$ with a faithful irreducible character
$\chi$ has cyclic center. We claim that if $G$ is a $p$-group of
nilpotency class $2$ and $\chi$ is a faithful irreducible character
then $N_w^\chi\geq 0$ if and only if $N_w(1)\geq N_w(z)$, where
$z\in Z(G)=Z$ has order $p$.  Indeed, it is well known that $\chi
(1)\chi=\eta^G$, where $\eta$ is a (faithful) linear character of
$Z(\chi)=Z$ (see for instance \cite[Theorem 2.31 and Problem 6.3]{Isaacs}). Then
$$N_w^\chi=\frac{1}{\chi (1)}(N_w,\eta^G)=\frac{1}{\chi (1)}({N_w}_{|Z},\eta)_{Z}.$$
If the order of $Z$ is $p^r$ and $z_1$ is a    generator   with
$z=z_1^{p^{r-1}}$  we have
\begin{equation}\label{sum}({N_w}_{|Z},\eta)_{Z}=\frac{1}{|Z|}\sum_{0< i\leq p^r} N_w(z_1^i)\varepsilon ^i=
\frac{1}{|Z|}\sum_{0\leq j\leq r}N_w(z_1^{p^j})(\sum_{ 0< i\leq
p^{r-j} \atop  (p , i)=1 } {(\varepsilon ^{p^j})}^i),\end{equation}
where $\varepsilon=\eta (z_1)$ is a primitive $p^r$-th root of unity
and we have used Lemma \ref{lem:coprime} for the second equality.
Notice that the innermost sum of (\ref{sum}) is the sum of all the
primitive $p^{r-j}$-th roots of unity, which is always zero except in
the cases $p^{r-j}=1$ or $p$, when it is $1$ or $-1$, respectively.
We conclude that
$$({N_w}_{|Z},\eta)_{Z}=\frac{1}{|Z|}(N_w(1)-N_w(z_1^{p^{r-1}})=\frac{1}{|Z|}(N_w(1)-N_w(z))$$
and our claim follows.

Now we prove the sufficiency part in the lemma. Let $\chi\in\Irr
(G)$, $K=\ker\chi$ and $G_1=G/K$. Of course we can suppose that
$\chi\not=1_G$ (because $N_w^{1_G}=|G|^{k-1}\geq 0$, $k$  is the
number of variables of $w$). By hypothesis $N_{w,G_1}(1)\geq
N_{w,G_1}(z)$, where $z$ is a central element of $G_1$ of order $p$.
We can view $\chi$   as a faithful character $\overline\chi$ of
$G_1$ and then our previous claim implies that
$N_{w,G_1}^{\overline\chi}\geq 0$. By Lemma \ref{coefswithkernel},
$N_w^\chi\geq 0$, which shows that $N_w$ is a character.

Conversely, suppose that $N_w$ is a character, that is, all the
numbers $N_w^\chi$ are non-negative, and consider an epimorphic
image $G_1=G/N$ with cyclic center and a central element $z\in G_1$
of order $p$. Then $G$ has an irreducible character $\chi$ with
kernel $N$ that is faithful when considered as a character
$\overline\chi$ of $G_1$. Since  $N_w^\chi\geq 0$, again by  Lemma
\ref{coefswithkernel}, $N_{w,G_1}^{\overline\chi}\geq 0$ and, by our
initial claim, the inequality  $N_{w,G_1}(1)\geq N_{w,G_1}(z)$
follows.
\end{proof}

\begin{theorem} Let $G$ be a $p$-group of nilpotency class $2$, $p$ odd,
 and $w$ a word.  Then $N_w$ is a character of $G$.
\end{theorem}

\begin{proof}
By the last result it suffices to show that if $G$ has cyclic center
$Z$ and $z\in Z$ has order $p$, $N_w(1)\geq N_w(z)$. Also we can
assume that $w$ has the form (\ref{noncommcanonicalword}) (if $w$ is as in (\ref{commcanonicalword}), skip the next two paragraphs).

If $Z^{p^{s_1}}\not=1$, we can write $z=z_1^{p^{s_1}}$ for some
$z_1\in Z$ and then the map $(g_1,g_2,\dots,g_k)\mapsto
(g_1z_1,g_2,\dots ,g_k)$ is a bijection between the sets of
solutions of $w=1$ and $w=z$, so $N_w(1)=N_w(z)$.

Now we suppose that $Z^{p^{s_1}}=1$ and notice that, since $G$ has nilpotency class $2$ and $p$ is
odd,
\begin{equation}\label{product}(xy)^{p^{s_1}}=x^{p^{s_1}}y^{p^{s_1}}[y,x]^{p^{s_1}\choose
2}=x^{p^{s_1}}y^{p^{s_1}}.\end{equation}

Therefore if we fix $g_2,g_4,\dots \in G$,
the map $(g_1,g_3,\dots )\mapsto w(g_1,g_2,\dots ,g_k)$ is a group
homomorphism $\varphi_{g_2,g_4,\dots}$. Obviously there is a
bijection between the kernel of this homomorphism and the set of
solutions of $w=1$ with $x_{2i}=g_{2i}$. As for the solutions of
$w=z$ with $x_{2i}=g_{2i}$, either this set is empty or else its
elements are in one-to-one correspondence with the elements in a
coset of the kernel of $\varphi_{g_2,g_4,\dots}$. In any case,
considering only solutions with $x_{2i}=g_{2i}$, the number of
solutions of $w=1$ is greater than or equal to the number of
solutions of $w=z$. Varying $g_2$, $g_4$,\dots, we conclude
$N_w(1)\geq N_ w(z)$, as desired.
\end{proof}

\section{The functions $N_{x^n}$ for $2$-groups of nilpotency class
$2$}

In this section we  study the functions $N_{x^n}$ for $2$-groups of
nilpotency class $2$ and characterize when this function is a
character. As we already pointed out in Section 3, the function
$N_{x^2\kern -1pt,\,Q_8}$ is not a character and in fact for each
$r\geq 1$ we can define a $2$-group $\QQ_{2^{3r}}$ of order
$2^{3r}$, which is the usual quaternion group $Q_8$ when $r=1$, such
that $N_{x^{2^r}\kern -4pt,\,\QQ_{2^{3r}}}$ is not a character. We
shall see that this group is in some sense involved in $G$ whenever
$N_{x^{2^r}\kern -4pt,\,G}$ is not a character. We shall also need
to introduce another family of groups, denoted $\DD_{2^{3r}}$, that, for
$r=1$, is the usual dihedral group of order $8$.

\begin{definition} Let $r\geq 1$. We define the \emph{quasi
dihedral} and \emph{quasi quaternion} group, $\DD_{2^{3r}}$ and
$\QQ_{2^{3r}}$, as
\begin{eqnarray*}
&& \DD_{2^{3r}}=\langle x,z_1,y\mid x^{2^r}=z_1^{2^r}=y^{2^r}=1,\, [x,z_1]=1,\, [x,y]=z_1,\, [z_1,y]=1\rangle,\\
&&\QQ_{2^{3r}}=\langle x,z_1,y\mid z_1^{2^r}=1,\,
x^{2^r}=z_1^{2^{r-1}}=y^{2^r},\, [x,z_1]=1,\, [x,y]=z_1,\,
[z_1,y]=1\rangle.
\end{eqnarray*}
\end{definition}

One can check that,  if $G=\DD_{2^{3r}}$ or $\QQ_{2^{3r}}$, $G$ has
order $2^{3r}$,  exponent $2^{r+1}$ and  $G'=Z(G)=\langle
z_1\rangle$ is cyclic of order $2^r$. Moreover, if $z=z_1^{2^{r-1}}$
is the central involution, in the (quasi) dihedral case
$N_{x^{2^r}}(1)=3\times 2^{3r-2}$ and $N_{x^{2^r}}(z)=  2^{3r-2}$,
whereas in the quaternion case the numbers are in the reverse order:
$N_{x^{2^r}}(1)=2^{3r-2}$ and $N_{x^{2^r}}(z)= 3\times 2^{3r-2}$
(and so $N_{x^{2^r}}$ is \emph{not} a character of $\QQ_{2^{3r}}$).

If $T$ and $H$ are $2$-groups with cyclic center and $|Z(T)|\leq
|Z(H)|$, we shall denote $T*H$ the central product of $T$ and $H$
with $Z(T)$ amalgamated with the corresponding subgroup of $Z(H)$.
Notice that if all the generators of $Z(T)$  are in the same orbit
under the action of the automorphism group of $T$ (or if a similar
situation holds in $H$), the group $T*H$ is unique up to isomorphism
and this is what happens if $T=\DD_{2^{3r}}$ or $\QQ_{2^{3r}}$.
Also, for a $p$-group $G$, $\Omega _r(G)$ is the subgroup generated
by the elements of order at most $p^r$.

\begin{lemma}\label{centralproduct} Let $G$ be a $2$-group of nilpotency class $2$ and cyclic center $Z$ of order $2^r$.
Suppose that $\Omega_{r+1}(G)'=Z$. Then $G=T*H$, where $T$ is
isomorphic to $\DD_{2^{3r}}$ or $\QQ_{2^{3r}}$.
\end{lemma}

\begin{proof}
Since $G$ has nilpotency class $2$, $\Omega_{r+1}(G)'$ is generated
by the commutators of elements of order at most $2^{r+1}$ and it is
cyclic, because it is contained in $Z$, which is cyclic, so  it is
generated by one of these commutators, say $[x,y]$. The orders of
$x$ and $y$ have to be $2^r$ or $2^{r+1}$ (because $[x,y]$ has order
$2^r$). If both have order $2^r$ it is clear that $T=\langle
x,y\rangle$ is isomorphic to $\DD_{2^{3r}}$ and, if both  have order
$2^{r+1}$, is isomorphic to $\QQ_{2^{3r}}$ (notice that
$G^{2^r}\subseteq Z$, so $x^{2^r}=y^{2^r}$). On the other hand, if
one is of order $2^r$ and the other of order $2^{r+1}$, their product
has  order $2^r$ and $T$   is isomorphic  to $\DD_{2^{3r}}$
again.

Now it suffices to check that $G=TC_G(T)$ (because obviously $T\cap
C_G(T)=Z(T)$ has order $2^r$, and so is the center of $G$). Indeed,
the conjugacy class of $x$ has order $|[x,G]|=|G'|=|Z|=p^r$ and
the same for $y$, so $$|G:C_G(T)|=|G:C_G(x)\cap C_G(y)|\leq
|G:C_G(x)| |G:C_G(y)|= 2^{2r}.$$ But
 $$|TC_G(T):C_G(T)|=|T:T\cap C_G(T)|=|T:Z| = 2^{2r},$$
so $G=TC_G(T)$, as claimed.
\end{proof}

One can check that, as it happens with the usual dihedral and
quaternion groups, $\DD_{2^{3r}}*\DD_{2^{3r}}$ and
$\QQ_{2^{3r}}*\QQ_{2^{3r}}$ are isomorphic. Using this result and
iterating Lemma \ref{centralproduct} we get  the following.

\begin{proposition}\label{centraldecomposition}Let $G$ be a $2$-group of nilpotency class $2$ and cyclic center $Z$ of order $2^r$.  Then $G$ is isomorphic to a group $\DD_{2^{3r}}*\stackrel{n}{\dots} *\DD_{2^{3r}}*H$ or $\DD_{2^{3r}}*\stackrel{n}{\dots} *\DD_{2^{3r}}*\QQ_{2^{3r}}*H$, $n\geq 0$, where $H$ has cyclic center of order $2^r$ and $\Omega_{r+1}(H)'$ is properly contained in the center of $H$.
\end{proposition}

Now suppose that $G=T*H$, where $T=\DD_{2^{3r}}$ or $\QQ_{2^{3r}}$
and $H$ is a $2$-group of nilpotency class $2$ and cyclic center
of order $2^r$. For any $g\in T$, $g^{2^r}=1$ or $z$, thus if $h\in
H$, $(gh)^{2^r}=1$ if and only if $g^{2^r}=h^{2^r}=1$ or $z$.
Similarly, $(gh)^{2^r}=z$ if and only if $g^{2^r}=1$ and $h^{2^r}=z$
or the other way round. This means that
\begin{eqnarray*}
&&N_{x^{2^r}\kern -4pt,\,G}(1)=(N_{x^{2^r}\kern
-4pt,\,T}(1)N_{x^{2^r}\kern -4pt,\,H}(1)+N_{x^{2^r}\kern
-4pt,\,T}(z)
N_{x^{2^r}\kern -4pt,\,H}(z))/2^r\\
&&N_{x^{2^r}\kern -4pt,\,G}(z)=(N_{x^{2^r}\kern
-4pt,\,T}(1)N_{x^{2^r}\kern -4pt,\,H}(z)+N_{x^{2^r}\kern
-4pt,\,T}(z)N_{x^{2^r}\kern -4pt,\,H}(1))/2^r,
\end{eqnarray*}
whence
\begin{eqnarray*}
&&N_{x^{2^r}\kern -4pt,\,G}(1)=2^{2r-2}(3N_{x^{2^r}\kern
-4pt,\,H}(1)+N_{x^{2^r}\kern -4pt,\,H}(z))
\text{\ \ or\ \ }  2^{2r-2}(N_{x^{2^r}\kern -4pt,\,H}(1)+3N_{x^{2^r}\kern -4pt,\,H}(z))\\
&&N_{x^{2^r}\kern -4pt,\,G}(z)=2^{2r-2}(3N_{x^{2^r}\kern
-4pt,\,H}(z)+N_{x^{2^r}\kern -4pt,\,H}(1)) \text{\ \ or\ \ }
2^{2r-2}(N_{x^{2^r}\kern -4pt,\,H}(z)+3N_{x^{2^r}\kern
-4pt,\,H}(1)),
\end{eqnarray*}
depending on whether $T=\DD_{2^{3r}}$ or $\QQ_{2^{3r}}$,
respectively. It follows that, in the former case, $N_{x^{2^r}\kern
-4pt,\,G}(1)\geq N_{x^{2^r}\kern -4pt,\,G}(z)$ if and only if
$N_{x^{2^r}\kern -4pt,\,H}(1)\geq N_{x^{2^r}\kern -4pt,\,H}(z)$ but,
in the latter case, this holds if and only if $N_{x^{2^r}\kern
-4pt,\,H}(1)\leq N_{x^{2^r}\kern -4pt,\,H}(z)$ (and the same
equivalences hold if inequalities are replaced by equalities). The
same happens if $T= \DD_{2^{3r}}*\stackrel{n}{\dots}
*\DD_{2^{3r}}$ or $\DD_{2^{3r}}*\stackrel{n}{\dots}
*\DD_{2^{3r}}*\QQ_{2^{3r}}$, $n\geq 0$.
The combination of this with Proposition \ref{centraldecomposition}
basically reduces our problem to the groups $G$ with
$\Omega_{r+1}(G)'$ properly contained in the center, which is the
situation considered in the next lemma.

\begin{lemma} \label{key} Let $G$ be a $2$-group of nilpotency class $2$ and cyclic center of order $2^r$ and let $z$ be the unique central involution. Suppose that $\Omega_{r+1}(G)'$ is properly contained in $Z=Z(G)$. Then $N_{x^{2^r}}(1)>N_{x^{2^r}}(z)$ if and only if $G$ has exponent $2^r$. Otherwise $N_{x^{2^r}}(1)=N_{x^{2^r}}(z)$.
\end{lemma}

\begin{proof}
Since $(G')^{2^r}\subseteq Z^{2^r}=1$,  raising to the  $2^{r+1}$-th
power is a group endomorphism of $G$ by (\ref{product}) and
$\Omega_{r+1}(G)=\{x\in G\mid x^{2^{r+1}}=1\}$. Moreover,
$\Omega_{r+1}(G)'$ is contained in $Z^2$,  so
$(\Omega_{r+1}(G)')^{2^{r-1}}=1$ and raising to the $2^r$-th power is
a group endomorphism of $\Omega_{r+1}(G)$ with kernel
$\Omega_r(G)=\{x\in G\mid x^{2^r}=1\}$. It is clear now that
$N_{x^{2^r}}(1)=|\Omega_r(G)|$ and
$N_{x^{2^r}}(z)=|\Omega_{r+1}(G)|-|\Omega_r(G)|$ (for any element
$x$ in $\Omega_{r+1}(G)\backslash  \Omega_r(G)$, $x^{2^r}$ is a central involution, so it is
$z$), so $N_{x^{2^r}}(1)>N_{x^{2^r}}(z)$ if and only if
$|\Omega_{r+1}(G):\Omega_r(G)|<2$, that is
$\Omega_{r+1}(G)=\Omega_r(G)=G$, i.~e., $G$ has exponent $2^r$.
Otherwise $1\not=\Omega_{r+1}(G)^{2^r}=\{x^{2^r}\mid
x\in\Omega_{r+1}(G)\}$, thus $z\in\Omega_{r+1}(G)^{2^r}$  is in the
image of the $2^r$-th power endomorphism of $\Omega_{r+1}(G)$ and
$N_{x^{2^r}}(z)=|\Omega_r(G)|=N_{x^{2^r}}(1)$.
\end{proof}

\begin{proposition}\label{prp:H} Let $G$ be a $2$-group of nilpotency class $2$, cyclic center of order $2^r$ and central involution $z$. Then $N_{x^{2^r}}(1)<N_{x^{2^r}}(z)$ if and only if $G$ is isomorphic to a group $\DD_{2^{3r}}*\stackrel{n}{\dots} *\DD_{2^{3r}}*\QQ_{2^{3r}}*H$, $n\geq 0$, where $H$ has cyclic center of order $2^r$ and exponent $2^r$.
\end{proposition}

\begin{proof} By Proposition \ref{centraldecomposition}, $G$ has two possible decompositions as a central product with one factor
 $H$ satisfying the hypotheses of Lemma  \ref{key}. If $\QQ_{2^{3r}}$ does not occur in the decomposition of $G$, we know that
 $N_{x^{2^r}}(1)<N_{x^{2^r}}(z)$ if and only if $N_{x^{2^r}\kern -4pt ,\,H}(1)<N_{x^{2^r}\kern -4pt ,\,H}(z)$, something that,
 according to Lemma \ref{key}, never happens. Therefore $\QQ_{2^{3r}}$ does occur in the decomposition of $G$. In this case,
 we know that $N_{x^{2^r}}(1)<N_{x^{2^r}}(z)$ if and only if $N_{x^{2^r}\kern -4pt ,\,H}(1)>N_{x^{2^r}\kern -4pt ,\,H}(z)$,
 which, by Lemma \ref{key}, is equivalent to $H$ having exponent $2^r$.
\end{proof}

It is not difficult to classify the groups $H$ in the previous
proposition.

\begin{lemma} Let $H$ be a $2$-group of nilpotency class $2$,
cyclic center of order $2^r$ and exponent $2^r$. Then $H\cong
\DD_{2^{3r_1}}*\dots *\DD_{2^{3r_n}}*C_{2^r}$ with $r_1\leq
\dots\leq r_n<r$.\end{lemma}

\begin{proof}
The result is trivially true if $H$ is abelian, so suppose
$H'=\langle [x,y]\rangle\not=1$ has order $2^s$. Since
$(xy)^{2^r}=1$, (\ref{product}) yields $s<r$. The elements $x^{2^s}$
and $y^{2^s}$ are central with orders at most $2^{r-s}$, so they lie
in $\langle z_1^{2^s}\rangle$, where $Z(H)=\langle z_1\rangle$, and, for
suitable $i$ and $j$, $xz_1^i$ and $yz_1^j$ have order exactly $2^s$. By
replacing $x$ and $y$ by these elements, we can suppose that
$T=\langle x,y\rangle\cong \DD_{2^{3s}}$. Arguing as in the last
part of the proof of Lemma \ref{centralproduct}, $H=TC_H(T)$ and
$T\cap C_H(T)=Z(T)$ is cyclic of order $2^s$. Since $Z(H)\leq
C_H(T)$, the hypotheses still hold in $C_H(T)$, so we can apply
induction.
\end{proof}

The last two results show that, with the hypotheses of Proposition
\ref{prp:H}, $N_{x^{2^r}}(1)<N_{x^{2^r}}(z)$ if and only if
$G\cong\DD_{2^{3r_1}}*\dots *\DD_{2^{3r_n}}*\QQ_{2^{3r}}$, $n\geq
0$, $r_1\leq\dots\leq r_n\leq r$. Notice simply that the cyclic
factor of $H$ is absorbed by $\QQ_{2^{3r}}$.

\begin{theorem} Let $G$ be a $2$-group of nilpotency class $2$. Then $N_{x^{2^r}}$ is a character of $G$
if and only if $G$ has no epimorphic image isomorphic to
$\DD_{2^{3r_1}}*\dots *\DD_{2^{3r_n}}*\QQ_{2^{3r}}$, $n\geq 0$,
$r_1\leq\dots\leq r_n\leq r$.
\end{theorem}

\begin{proof}
If $G$ has an epimorphic image $G_1$ of the indicated type, then by
the last remark, $N_{x^{2^r}\kern -4pt ,\,G_1}(1)<N_{x^{2^r}\kern
-4pt ,\,G_1}(z)$ ($z$ is the central involution of $G_1$) and, by
Lemma \ref{basic}, $N_{x^{2^r}}$ is not a character of $G$.
Conversely, if $N_{x^{2^r}}$ is not a character of $G$, by the same
lemma, $G$ has an epimorphic image $G_1$ with cyclic center such
that $N_{x^{2^r}\kern -4pt ,\,G_1}(1)<N_{x^{2^r}\kern -4pt
,\,G_1}(z)$. We claim that the center $Z$ of $G_1$ has order $2^r$
(and then, again by the last remark, $G_1$ is the desired epimorphic
image). The map $x\mapsto x^{2^r}$ cannot be a group endomorphism of
$G_1$ (this would immediately imply that either $N_{x^{2^r}\kern
-4pt ,\,G_1}(z)=0$ or else $N_{x^{2^r}\kern -4pt
,\,G_1}(1)=N_{x^{2^r}\kern -4pt ,\,G_1}(z)$), hence
$(G_1')^{2^{r-1}}\not=1$ and $Z^{2^{r-1}}\not=1$, that is $|Z|\geq
2^r$ (here $Z$ is the center of $G_1$). If $|Z|>2^r$, $z=z_1^{2^r}$
for some $z_1\in Z$ and then $N_{x^{2^r}\kern -4pt
,\,G_1}(1)=N_{x^{2^r}\kern -4pt ,\,G_1}(z)$ because $x\mapsto xz_1$
maps bijectively the solutions of $x^{2^r}=1$ to the solutions of
$x^{2^r}=z$. Thus $|Z|=2^r$ and the proof is complete.
\end{proof}

\section{$p$-groups with central Frattini subgroup}

In this section we consider a $d$-generated $p$-group $G$ such that
$\Phi(G)\leq Z(G)$, that is, $G$ has nilpotency class $2$ and
elementary  abelian derived subgroup. We shall show that for the words
$w_k=[x_1,x_2]\dots [x_{2k-1},x_{2k}]$, with $k\geq d_0=\lfloor d/2\rfloor$,
\begin{equation}\label{fiber}N_{w_k}(g)\geq |G|^{2k-1}\text{\ \  for all }g\in G_{w_k}.\end{equation}
For $k=1$, (\ref{fiber}) can be easily proved for any $p$-group of
nilpotency class $2$. Indeed, if $g=[x,y]$, we can multiply $x$ by
any element commuting with $y$ and $y$ by any central element, thus
$N_{w_1}(g)\geq |C_G(y)||Z(G)|=|G||Z(G):[y,G]|\geq |G|$.

If $V=G/\Phi(G)$, viewed as a vector space over $\F_p$, there is a
natural surjective linear map $\pi$ from $\bigwedge^2=\bigwedge^2V$,
the exterior square of $V$, onto $G'$ given by $\overline
x\wedge\overline y\mapsto [x,y]$.
 For a fixed $\omega\in \bigwedge^2$ and $k\geq 0$, it is then natural to consider the number
 $N_{w_k}(\omega)=N_{w_k,V}(\omega)$ of solutions in $V^{(2k)}$ of the equation
\begin{equation}\label{vectorword}w_k(x_1,\dots,x_{2k})=x_1\wedge x_2+\cdots +x_{2k-1}\wedge x_{2k}=\omega.
\end{equation}
The set of values of $w_k$, that is, the set
$\{w_k(v_1,\dots,v_{2k})\mid v_i\in V\}\subseteq\bigwedge^2$ will be
denoted  simply $\bigwedge^2_{w_k}$.

Similarly as in Section 2, if we fix a basis $\{e_1,\dots,e_d\}$ of
$V$ there is a one-to-one correspondence between $\bigwedge^2$ and
$\AA_d$, the set of $d\times d$ antisymmetric matrices over the
field $\F_p$, given by $\sum_{1\leq i<j\leq d}a_{ij}(e_i\wedge
e_j)\mapsto A$, where $A\in\AA_d$ has entries $a_{ij}$ for $1\leq
i<j\leq d$. Then solving (\ref{vectorword}) amounts to solving the
matrix equation $X^tJ_kX=A$, where $X$ represents a $2k\times d$
matrix, $J_k$ is the $2k\times 2k$ block diagonal matrix with
repeated diagonal block
$\begin{pmatrix}\phantom{-}0&1\\-1&0\end{pmatrix}$ and $A\in\AA_d$.
Moreover, the number of solutions only depends on the rank of $A$,
so there is no loss to assume that
$A=\begin{pmatrix}J_r&0\\0&0\end{pmatrix}$ with $r\leq k$ (otherwise
there are no solutions). This number, which is denoted
$N(K_{2k},K_{d,2r})$ in  \cite{WeiZhang}, was originally considered
by Carlitz \cite{Carlitz} and later on by other authors (see
\cite{WeiZhang}).  If $\omega\in
\bigwedge^2_{w_r}\backslash\bigwedge^2_{w_{r-1}}$, the corresponding
antisymmetric matrix $A$ has rank $2r$, so the number
$N(K_{2k},K_{d,2r})$ is, in our notation, $N_{w_k}(\omega)$. An
explicit formula for this number is given in \cite[Theorems 3,4,
5]{WeiZhang} that, for $r\leq k\leq d_0=\lfloor d/2\rfloor$, can be
written as
\begin{equation}\label{Nformula}
N(K_{2k},K_{d,2r})=p^{r(2k-r)}\prod_{i=k-r+1}^k(p^{2i}-1)\sum_{j=0}^{k-r}p^{j\choose 2}\left[\begin{matrix}d-2r\\ j\end{matrix}\right]_p\prod_{i=k-r-j+1}^{k-r}(p^{2i}-1),\end{equation}
where  $\left[\begin{matrix}n\\ j\end{matrix}\right]_p=
(p^n-1)\dots (p^{n-j+1}-1)/(p^j-1)\dots (p-1)$ for $j>0$ and
$\left[\begin{matrix}n\\ 0\end{matrix}\right]_p=1$. If we just consider in (\ref{Nformula}) the summand corresponding to $j=k-r$ we get the inequality
 \begin{equation}\label{Ninequality}
N(K_{2k},K_{d,2r})\geq p^{r(2k-r)}\prod_{i=1}^k(p^{2i}-1)p^{k-r\choose 2}\left[\begin{matrix}d-2r\\ k-r\end{matrix}\right]_p.\end{equation}
But $\prod_{i=1}^k(p^{2i}-1)>\prod_{i=1}^k p^{2i-1}=p^{k^2}$, and $\left[\begin{matrix}n\\ j\end{matrix}\right]_p>(\prod_{i=n-j+1}^{n}p^{i-1})/(\prod_{i=1}^jp^i)=p^{nj}$, \smash{therefore}
$$N(K_{2k},K_{d,2r})> p^{r(2k-r)+k^2+{k-r\choose 2}+(d-2r)(k-r)}.$$
The quadratic function of $r$ in the exponent of $p$ in this formula is decreasing in the interval $0\leq r\leq k$, so we conclude that
\begin{equation}\label{Nestimation}N_{w_k}(\omega)=N(K_{2k},K_{d,2r})> p^{2k^2}  \text{ for any }\omega\in\textstyle\bigwedge\nolimits^2_{w_k}.\end{equation}

Since the rank of a $d\times d$ antisymmetric matrix is at most
$2d_0$ we have that $\bigwedge^2_{w_k}=\bigwedge^2$ for $k\geq d_0$.
But for any $k$, $\pi$ maps $\bigwedge^2_{w_k}$ onto $G_{w_k}$, thus
$G_{w_k}=G'$ for $k\geq d_0$. Now it is easy to show that if
(\ref{fiber})  holds for $w=w_{d_0}$, it also holds for  $w=w_k$ for
any $k\geq d_0$. Indeed, if $k\geq d_0$, we have
$N_{w_k}=N_{w_{d_0}}*N_{w_{k-d_0}}$ and, since we are assuming that
$N_{w_{d_0}}(x)\geq |G|^{2d_0-1}$ for any $x\in G'$, if $g\in G'$,
\begin{eqnarray*}N_{w_k}(g)&=&\sum_{y\in G_{w_{k-d_0}}}N_{w_{d_0}}(gy^{-1})N_{w_{k-d_0}}(y)\geq
|G|^{2d_0-1}\sum_{y\in
G_{w_{k-d_0}}}N_{w_{k-d_0}}(y)\\&=&|G|^{2d_0-1}|G|^{2(k-d_0)}=|G|^{2k-1}.\end{eqnarray*}
So in order to prove (\ref{fiber}) for $w_k$ we can always assume
that $1\leq k\leq d_0$.

It is clear that if $g\in G'$ and $\pi(\omega)=g$, the solutions of
$w_k(x_1,\dots,x_{2k})=\omega$ in $V$ can be lifted to solutions of
$w_k(x_1,\dots,x_{2k})=g$ in $G$ and of course, all solutions of the
equation in $G$ occur in this way, so
$$N_{w_k}(g)=|\Phi(G)|^{2k}\sum_{\omega\in\pi^{-1}(g)}N_{w_k,V}(\omega)$$
and (\ref{fiber}) can be written now as
\begin{equation}\label{Amit}|\Phi(G)|\sum_{\omega\in\pi^{-1}(g)}N_{w_k,V}(\omega)\geq |G:\Phi(G)|^{2k-1}=p^{d(2k-1)}\end{equation}
for $g\in G_{w_k}$. Obviously only the $\omega$'s in $\pi^{-1}(g)\cap\bigwedge^2_{w_k}$ contribute to this sum and for them we can use the estimation (\ref{Nestimation}). Thus, since $|G'|=p^{d(d-1)/2}/|\ker\pi|\leq |\Phi (G)|$, the inequality
\begin{equation}\label{Amit2}p^{\frac{d(d-1)}{2}+2k^2-d(2k-1)}|\pi^{-1}(g)\cap\textstyle\bigwedge\nolimits^2_{w_k}|\geq |\ker\pi|\end{equation}
implies (\ref{Amit}). If $k=d_0$, $|\pi^{-1}(g)\cap\textstyle\bigwedge\nolimits^2_{w_k}|=|\pi^{-1}(g)|=|\ker\pi|$ and (\ref{Amit2}) holds because
the exponent of $p$ is positive (it is positive for any $k$). We conclude the following result.

\begin{proposition} Let $G$ be a $d$-generated $p$-group with $\Phi(G)\leq Z(G)$. Then for any $k\geq \lfloor d/2\rfloor$ and $g\in G'$, $N_{w_k}(g)\geq |G|^{2k-1}$.
\end{proposition}

Another situation in which $|\pi^{-1}(g)\cap\textstyle\bigwedge\nolimits^2_{w_k}|=|\ker\pi|$ is when $\pi$ is an isomorphism, so we also have the following result.

\begin{proposition} Let $G$ be a $d$-generated $p$-group with $\Phi(G)\leq Z(G)$ and $|G'|=p^{d(d-1)/2}$. Then for any $k\geq 1$ and  $g\in G_{w_k}$, $N_{w_k}(g)\geq |G|^{2k-1}$.
\end{proposition}

Notice that the last proposition applies in particular to the free
$p$-groups of nilpotency class $2$ and exponent $p$ and, in this
case the inequality $N_w(g)\geq |G|^{k-1}$,  $g\in G_w$, is in fact
true for any word $w\in F_k$. This is clear if $w\in F_k'$ and
otherwise we can suppose that $s_1=0$ in
(\ref{noncommcanonicalword}). But in this case all equations $w=g$
have the same number of solutions, namely, $|G|^{2k-1}$.


\begin{thebibliography}{00}

\bibitem{AmitVishne}
A.~Amit, U.~Vishne, Characters and solutions to equations in finite groups,
J. Algebra Appl. {\bf 10}, no. 4, (2011) 675--686.

\bibitem{Carlitz}
L.~Carlitz,  Representations by skew forms in a finite field, Arch. Math. {\bf 5} (1954) 19--31.

\bibitem{DasNath}
A.~K.~Das, R.~K.~Nath, On solutions of a class of equations in a finite group, Comm. Algebra {\bf 37}, no. 11, (2009) 3904--3911.

\bibitem{Isaacs}
I.~M.~Isaacs, \emph{Character Theory of Finite Groups}, Dover
Publications INC., New York, 1994.

\bibitem{Khukhro}
E.~Khukhro, \emph{$p$-automorphisms of finite $p$-groups}, Cambridge University Press, Cambridge, 1998.

\bibitem{Levy}
M.~Levy, On the probability of satisfying a word in
nilpotent groups of class 2, arXiv:1101.4286 (2011).

\bibitem{Lubotzky}
A.~Lubotzky, Images of word maps in finite simple groups, arXiv:1211.6575 (2012).

\bibitem{Newman}
M.~Newman, \emph{Integral matrices}, Academic Press, New York, 1972.

\bibitem{NikolovSegal}
N.~Nikolov, D.~Segal, A characterization of finite soluble groups, Bull. Lond. Math. Soc. {\bf 39} (2007) 209--213.

\bibitem{Parzanchevski}
O.~Parzanchevski, G.~Schul, On the Fourier expansion of word maps, Bull. Lond. Math. Soc. {\bf 46} (2014) 91--102.

\bibitem{Solomon}
L.~Solomon, The solution of equations in groups, Arch. Math. {\bf
20} (1969) 241--247

\bibitem{Strunkov}
S.~P.~Strunkov, On the theory of equations on finite groups,  Izv.
Math. {\bf 59}, no. 6, (1995) 1273--1282.

\bibitem{WeiZhang}
J.~Wei, Y.~Zhang, The number of solutions to the alternate matrix equation over a finite field and a $q$-identity,
J. Statist. Plann. Inference {\bf 94},  no. 2, (2001) 349--358.


\end{thebibliography}
\end{document}